\documentclass[12pt, reqno]{amsart}
\usepackage{amssymb, mathrsfs}
\usepackage{latexsym}
\usepackage{amsmath}
\usepackage{amsthm}
\usepackage{amsfonts}
\usepackage{dsfont, upgreek}
\usepackage{mathtools}
\usepackage{epsfig}
\usepackage{cite}
\usepackage{amscd}
\usepackage{graphicx}
\usepackage{tikz-cd}
\usepackage{enumitem}
\setlist[enumerate]{itemsep=0.5ex}

\usepackage{stmaryrd}

\usepackage{color}
\usepackage{tikz}
\usetikzlibrary{patterns}

\usepackage[left=1.4in,right=1.4in,top=1.2in,bottom=1.2in]{geometry}

\usepackage{tikz}
\usetikzlibrary{decorations.markings}
\usetikzlibrary{arrows.meta}

\usepackage{extarrows}

\usepackage{adjustbox}
\usepackage{pgfplots}
\usepgfplotslibrary{colormaps}
\usepackage{slashed}

\usepackage[colorlinks=true,linkcolor=blue,citecolor=blue]{hyperref}

\usepackage[colorinlistoftodos,prependcaption,textsize=tiny]{todonotes}  

\usepackage[all,cmtip]{xy}
\theoremstyle{plain}
\newtheorem{theorem}{Theorem}[section]

\newtheorem{lemma}[theorem]{Lemma}

\newtheorem{question}[theorem]{Question}

\theoremstyle{definition}

\newtheorem*{claim*}{Claim}

\theoremstyle{remark}

\numberwithin{equation}{section}

\newcommand{\R}{\mathbb{R}}

\newcommand{\interior}[1]{%
	{\kern0pt#1}^{\mathrm{\,o}}%
}

\makeatletter
\let\save@mathaccent\mathaccent
\newcommand*\if@single[3]{%
	\setbox0\hbox{${\mathaccent"0362{#1}}^H$}%
	\setbox2\hbox{${\mathaccent"0362{\kern0pt#1}}^H$}%
	\ifdim\ht0=\ht2 #3\else #2\fi
}
\newcommand*\rel@kern[1]{\kern#1\dimexpr\macc@kerna}
\newcommand*\wideaccent[2]{\@ifnextchar^{{\wide@accent{#1}{#2}{0}}}{\wide@accent{#1}{#2}{1}}}
\newcommand*\wide@accent[3]{\if@single{#2}{\wide@accent@{#1}{#2}{#3}{1}}{\wide@accent@{#1}{#2}{#3}{2}}}
\newcommand*\wide@accent@[4]{%
	\begingroup
	\def\mathaccent##1##2{%
		\let\mathaccent\save@mathaccent
		\if#42 \let\macc@nucleus\first@char \fi
		\setbox\z@\hbox{$\macc@style{\macc@nucleus}_{}$}%
		\setbox\tw@\hbox{$\macc@style{\macc@nucleus}{}_{}$}%
		\dimen@\wd\tw@
		\advance\dimen@-\wd\z@
		\divide\dimen@ 3
		\@tempdima\wd\tw@
		\advance\@tempdima-\scriptspace
		\divide\@tempdima 10
		\advance\dimen@-\@tempdima
		\ifdim\dimen@>\z@ \dimen@0pt\fi
		\rel@kern{0.6}\kern-\dimen@
		\if#41
		#1{\rel@kern{-0.6}\kern\dimen@\macc@nucleus\rel@kern{0.4}\kern\dimen@}%
		\advance\dimen@0.4\dimexpr\macc@kerna
		\let\final@kern#3%
		\ifdim\dimen@<\z@ \let\final@kern1\fi
		\if\final@kern1 \kern-\dimen@\fi
		\else
		#1{\rel@kern{-0.6}\kern\dimen@#2}%
		\fi
	}%
	\macc@depth\@ne
	\let\math@bgroup\@empty \let\math@egroup\macc@set@skewchar
	\mathsurround\z@ \frozen@everymath{\mathgroup\macc@group\relax}%
	\macc@set@skewchar\relax
	\let\mathaccentV\macc@nested@a
	\if#41
	\macc@nested@a\relax111{#2}%
	\else
	\def\gobble@till@marker##1\endmarker{}%
	\futurelet\first@char\gobble@till@marker#2\endmarker
	\ifcat\noexpand\first@char A\else
	\def\first@char{}%
	\fi
	\macc@nested@a\relax111{\first@char}%
	\fi
	\endgroup
}
\makeatother

\makeatletter
\newcommand*{\transpose}{%
	{\mathpalette\@transpose{}}%
}
\newcommand*{\@transpose}[2]{%
	\raisebox{\depth}{$\m@th#1\intercal$}%
}
\makeatother

\begin{document}
\title[Infinitely many convex polytopes with dihedral angle upper bound]{Infinitely many convex polytopes with dihedral angles bounded above by $\pi/2+\delta$}
%
%
\author{Zhenkun Li}
\address[Zhenkun Li]{Institute of Mathematics, Chinese Academy of Sciences}
\email{zhenkun@amss.ac.cn}
\thanks{}
\author{Jinmin Wang}
\address[Jinmin Wang]{Institute of Mathematics, Chinese Academy of Sciences}
\email{jinmin@amss.ac.cn}
\thanks{}
\author{Jing Yu}
\address[Jing Yu]{Shanghai Center for Mathematical Sciences, Fudan University}
\email{jyu@fudan.edu.cn}
\thanks{}

\begin{abstract}
For every $\delta>0$, we construct infinitely many pairwise combinatorially distinct convex three-dimensional polytopes whose dihedral angles are all less than $\pi/2+\delta$.
\end{abstract}
\maketitle
\section{Introduction}
Throughout this paper, dihedral angles are understood to be interior dihedral angles. As reported by Adiprasito \cite{Adiprasito2022}, Gromov asked the following finiteness question:
\begin{question}\label{ques: maximal epsilon_d}
For a fixed dimension $d\geq 3$ and $\epsilon>0$, are there only finitely many combinatorial types of convex $d$-polytopes whose dihedral angles are all less than $\pi-\epsilon$?
\end{question}

Classical results imply finiteness for $\epsilon\geq\pi/2$; see \cite[Theorem 2]{Coxeter34}. Gromov discusses infinite families for small $\epsilon>0$ in \cite[Section 1.1]{Gromov24}, attributing the initial observation to Adiprasito. The corresponding question about the threshold for finiteness appears as Problem 18 in \cite[Page 196]{AdiprasitoGoaocPatakova2024}, attributed to Adiprasito.

We prove that $\pi/2$ is the optimal dihedral angle upper bound for every $d\geq3$. More precisely, we establish the following results.

\begin{theorem}\label{thm:example}
	For every $\delta>0$ and every positive integer $N$, there exists a convex three-dimensional polytope with at least $N$ two-dimensional faces whose dihedral angles are all less than $\pi/2+\delta$.
\end{theorem}

The following finiteness theorem (which directly answers Question \ref{ques: maximal epsilon_d} with the optimal $\epsilon$) follows from Theorem \ref{thm:example} and \cite[Theorem 2]{Coxeter34}.
\begin{theorem}\label{thm:finiteness}
	For any $d\geq 3$ and $\beta>0$, let $\mathcal P_{d,\beta}$ be the collection of convex $d$-dimensional polytopes with dihedral angles at most $\beta$, modulo combinatorial equivalence. Then
	\begin{enumerate}
		\item $|\mathcal P_{d,\pi/2}|<\infty$, and
		\item $|\mathcal P_{d,\pi/2+\delta}|=\infty$ for any $\delta>0$.
	\end{enumerate}
\end{theorem}
\begin{proof}
	\noindent\textup{(1)} By \cite[Theorem 2]{Coxeter34}, any member of $\mathcal P_{d,\pi/2}$ is a direct product of simplices. Thus $\mathcal P_{d,\pi/2}$ is finite. In particular, $|\mathcal P_{d,\pi/2}|$ equals the partition number of $d$, number of ways of writing $d$ as a sum of positive integers.
    
  \par\medskip\noindent\textup{(2)} For any $\delta>0$, $|\mathcal P_{d,\pi/2+\delta}|=\infty$ is witnessed by the direct product of the examples in Theorem \ref{thm:example} with $[0,1]^{d-3}$. 
\end{proof}

We prove Theorem \ref{thm:example} in Section \ref{sec:2}.

\vspace{0.2in}
\noindent{\bf AI disclosure}.
The authors used GPT-6 Astra to help develop the proof. The authors checked and rewrote the proof, and take full responsibility for its final form and correctness.

\vspace{0.2in}
\noindent{\bf Acknowledgement}.
Jinmin Wang is partially supported by NSFC 12501169.

\section{The construction in dimension three}\label{sec:2}
In this section, we will construct infinitely many combinatorially distinct polytopes with uniformly bounded dihedral angles. 

Below we will describe a special family $\{\vec{v}_i\}$ of normal vectors that supplies our construction.

Fix $R>1$ a large transcendental real number and $M>0$ a large integer. Define
\begin{equation}\label{eq: theta}
	\theta=\arccos\left(-1/R\right) \in \left(\frac{\pi}{2},\pi\right).
\end{equation}
For $j=0,1,\cdots,M-1$, define
\begin{equation}\label{eq: v_j}
	\vec{v}_j=(\cos(j\theta),\sin(j\theta),R^{j-M}), \quad \vec{n}_j = \frac{\vec{v}_j}{|\vec{v}_j|}.
\end{equation}

For suitable positive real values $h_j$ for $j=0,\cdots,M-1$, we construct a polytope $P_k$, for $k = 2,\cdots, M-1$, exactly enclosed bottom face $B=P_k\cap\{x_3=0\}$, with unit outer normal vector $\textbf{n}=(0,0,-1)$, and faces given by
$$H_j=\{\vec x\in\R^3:\langle\vec x,\vec v_j\rangle=h_j\}.$$
More precisely,
\begin{equation}\label{eq:Pk}
	P_k=\{\vec{x}\in \R^3:\langle \vec{x},\textbf{n}\rangle\leq 0\}\cap \Big(\bigcap_{j=0}^k\{\vec{x}\in\R^3:\langle \vec{v}_j,\vec{x}\rangle\leq h_j\}\Big).
\end{equation}
For a fixed $P_k$, write
\[
F_j=P_k\cap H_j\quad(0\le j\le k).
\]
suppressing the dependence of these faces on $k$. We will ensure that $P_k$ is bounded and that $B,F_0,\ldots,F_k$ are exactly its faces, with outward unit normals $\textbf{n},\vec n_0,\ldots,\vec n_k$, respectively.

 The basic idea is an induction: we start from a tetrahedron $P_2$ enclosed by $B$ and $F_0,F_1,F_2$. At each induction step, we cut $P_{k-1}$ with the plane $H_k$ near its ``top vertex'' incident with $F_{k-1},F_{k-2},F_{k-3}$, and get $F_k$.  We then control all newly created dihedral angles and show that they all stay bounded. To proceed, we need some linear algebra preparations.
\begin{lemma}\label{lemma:linearlyIndep} Let $\textbf{n} = (0,0,-1)$ and $\vec{v}_j$ be defined as in \eqref{eq: v_j}. Then the following statements hold.
	\begin{enumerate}
		\item For distinct nonnegative integers $i$ and $j$, the vectors $\textbf{n},\vec{v}_i,\vec{v}_j$ are linearly independent.
		\item For every integer $j \ge 0$, the vectors $\vec{v}_{j},\vec{v}_{j+1},\vec{v}_{j+2}$ are linearly independent.
	\end{enumerate}
\end{lemma}
\begin{proof}
\noindent\textup{(1)} A direct computation shows that $\det[\textbf{n},\vec{v}_i,\vec{v}_j]=\sin((i-j)\theta)$. This is a non-zero transcendental number by \eqref{eq: theta}.

\noindent\textup{(2)} We also have direct computation
		\begin{align*}
	\det[\vec{v}_{j},\vec{v}_{j+1},\vec{v}_{j+2}]=&~R^{j-M}\left(\sin\theta-R\sin (2\theta)+R^{2}\sin\theta\right)
		\\=&~R^{j-M}\left(R^2+3\right)\sin\theta,
		\end{align*}
		which is non-zero.
\end{proof}

A core ingredient of the construction is the following observation of a positive linear combination.
\begin{lemma}\label{lemma:positiveLinearCombination}
	For every integer $j\geq 0$,
	\[
	\vec{v}_{j+3}=R\vec{v}_{j}+\vec{v}_{j+1}+\left(R-\frac 2 R\right)\vec{v}_{j+2}.
	\]
\end{lemma}
\begin{proof}
	By \eqref{eq: theta}, $\cos\theta=-1/R$, thus we have
	\begin{equation*}
\cos((j+2)\theta)+\cos(j\theta)=2\cos((j+1)\theta)\cos\theta=-\frac 2 R\cos((j+1)\theta),
	\end{equation*}
	\begin{equation*}
\sin((j+2)\theta)+\sin(j\theta)=2\sin((j+1)\theta)\cos\theta=-\frac 2 R\sin((j+1)\theta).
	\end{equation*}
    Recall that 
    \[\vec{v}_j=(\cos(j\theta),\sin(j\theta),R^{j-M}).\]
	It follows that
	\[
		\vec{v}_{j+2}+\vec{v}_{j}+\frac 2 R \vec{v}_{j+1}=(0,0, R^{j-M}\left(R^2+3\right)).
	\]
	Therefore,
	\[
		R\Big(\vec{v}_{j+2}+\vec{v}_{j}+\frac 2 R \vec{v}_{j+1}\Big)= \vec{v}_{j+3}+\vec{v}_{j+1}+\frac 2 R \vec{v}_{j+2}.
	\]
Rearranging gives the required recurrence. 
\end{proof}

The following lemma will provide the desired dihedral angle bounds.
\begin{lemma}\label{lem: dihedral}
	For integers $0 \leq i < j \leq M-1$ with $j - i \le 3$, we have
	\[
		\langle \vec n_i,\vec n_j\rangle=\frac{\langle \vec{v}_i,\vec{v}_j\rangle}{|\vec{v}_i| |\vec{v}_j|}< \frac{4}{R}.
	\]
\end{lemma}
\begin{proof}
	A direct computation shows
	\[
		\frac{\langle \vec{v}_i,\vec{v}_j\rangle}{|\vec{v}_i| |\vec{v}_j|}=\frac{\cos((j-i)\theta)+R^{i+j-2M}}{\sqrt{(1+R^{2i-2M})(1+R^{2j-2M})}}.
	\]
	As $\cos\theta=-1/R$ and $1\leq j - i \leq 3$, the term $\cos((j-i)\theta)$ only takes one of the three values:
	\[
		\cos\theta=-\frac 1 R,~\cos(2\theta)=\frac{2}{R^2}-1,~\cos(3\theta)=\frac 3 R-\frac{4}{R^3}.
	\]
	In particular, we have
	\[
		\cos((j-i)\theta)<\frac 3 R.
	\]
	Furthermore, as $i,j\leq M-1$, we have
	\[
		\frac{R^{i+j-2M}}{\sqrt{(1+R^{2i-2M})(1+R^{2j-2M})}}=
		\frac{R^{i+j}}{\sqrt{(R^{2M}+R^{2i})(R^{2M}+R^{2j})}}< \frac{1}{R}.
	\]
	As a consequence,
	\[
		\frac{\langle \vec{v}_i,\vec{v}_j\rangle}{|\vec{v}_i| |\vec{v}_j|} < \frac{3}{R} + \frac{1}{R} = \frac{4}{R}.
	\]
    This concludes the proof.
\end{proof}

	\begin{proof}[Proof of Theorem \ref{thm:example}]
We now present the construction. For two faces $F,G$ sharing an edge, write $\angle(F,G)$ for their interior dihedral angle. 

	\noindent{\bf Base case}. When $k=2$, taking $h_0 = R^{-M}$, $h_1 = R^{1-M}$, and $h_2 = R^{2-M}$, we introduce $P_2$ as in \eqref{eq:Pk}. Thus $P_2$ is a tetrahedron, with a unique vertex above the $xy$-plane. Its bottom face is $B$ and side faces are $F_0,F_1,F_2$. The side dihedral angles are given by 
    $$\cos(\pi-\angle(F_i,F_j))=\langle\vec n_i,\vec n_j\rangle$$
    for $0\leq i<j\leq 2$.    By Lemma \ref{lem: dihedral},
	\[
		\begin{aligned}
			\angle(F_i,F_j) &= \arccos(-\langle\vec n_i,\vec n_j\rangle)< \arccos(-\frac{4}{R})= \frac{\pi}{2} + \arcsin(\frac{4}{R}).
		\end{aligned}
	\]
    For bottom dihedral angles, we have
    $$\cos(\pi-\angle(B,F_j))=\langle\textbf{n},\vec n_j\rangle$$
    for $j=0,1,2$. A direct computation shows
	\[
		\langle\textbf{n},\vec n_j\rangle  = -\frac{R^{j-M}}{\sqrt{1+R^{2j-2M}}} < 0.
	\]
	Hence $\angle(B,F_j)$ are all acute.

\vspace{0.1in}
\noindent{\bf Induction hypotheses}. 
We will construct inductively the polytopes $P_k$ for $k=2,3,\ldots,M-1$ such that the following properties hold.
\begin{enumerate}
	\item\label{ih:facets} $P_{k}$ has exactly $(k+2)$ faces,namely
\[
  B=P_k\cap\{x_3=0\}\quad\text{and}\quad
  F_j=P_k\cap H_j\quad(0\le j\le k).
\]
	\item\label{ih:simple} The degree of each vertex of $P_{k}$ is three.
	\item\label{ih:bottom} For $0\le j\le k$, $B$ and $F_j$ share an edge if and only if $j=0,1,2$.
	\item \label{ih:adjacency} For $0\le i<j\le k$, $F_i$ and $F_j$ share an edge if and only if $1\leq |i-j|\leq 3$.
	\item \label{ih:angles} All dihedral angles of $P_{k}$ are smaller than $(\pi/2+\arcsin(4R^{-1}))$.
    \item \label{ih:counts} $P_k$ has $2k$ vertices and $3k$ edges.
\end{enumerate}
The base case verifies all these properties for $P_2$.

\vspace{0.1in}
\noindent{\bf Induction step}. 
Suppose that $3\le k\le M-1$ and that $P_{k-1}$ has been constructed.
We now construct $P_k$ out of $P_{k-1}$. Let $p$ be the vertex of $P_{k-1}$ meeting three faces $F_{k-1},F_{k-2},F_{k-3}$, guranteed by Lemma \ref{lemma:linearlyIndep}. We cut the vertex $p$ away by a plane with normal vector $\vec v_k$.

\vspace{0.1in}
\noindent{\bf Claim}. 
The linear functional ${x}\mapsto\langle\vec{v}_k,{x}\rangle$ attains its maximum on $P_{k-1}$ uniquely at $p$.

Indeed, by Lemma \ref{lemma:positiveLinearCombination}, $\vec{v}_k$ is a positive linear combination of $\vec{v}_{k-1},\vec{v}_{k-2},\vec{v}_{k-3}$. Therefore, for any ${x}\in P_{k-1}$,
	\begin{align*}
		\langle \vec{v}_k,{x}\rangle\leq&~ 
        \left(R-\frac 2 R\right)\langle \vec{v}_{k-1},{x}\rangle+\langle \vec{v}_{k-2},{x}\rangle+R\langle \vec{v}_{k-3},{x}\rangle\\
		\leq&~\left(R-\frac 2 R\right)h_{k-1}+h_{k-2}+Rh_{k-3}.
	\end{align*}
	Equality holds if and only if 
	$$\langle \vec{v}_{k-1},{x}\rangle=h_{k-1},~\langle \vec{v}_{k-2},{x}\rangle=h_{k-2},~\langle \vec{v}_{k-3},{x}\rangle=h_{k-3},$$
	that is, ${x}=p$. This finishes the proof of the claim.

Set
\begin{equation}\label{eq:cut}
  h_k=\frac 1 2\big(\langle\vec v_{k},p\rangle+\max_{w\ne p}|\langle\vec v_{k},w\rangle|\big).
\end{equation}
where the maximum is taken for all vertices $w$ of $P_{k-1}$ that is not $p$. Consequently, the plane $H_k$ strictly separates $p$ from every other vertex of $P_{k-1}$. We define
\[ P_k=P_{k-1}\cap\{x\in\R^3:\langle\vec{v}_k,x\rangle\le h_k\},\]
in agreement with~\eqref{eq:Pk}.

Every old edge survives, possibly shortened, and $F_k$ is adjacent precisely to $F_{k-1},F_{k-2},F_{k-3}$. All old faces remain faces, and all old face adjacencies are preserved. Since $p\notin B$, the bottom face is unchanged. This verifies~\textup{(\ref{ih:facets})}-\textup{(\ref{ih:adjacency})}.

The truncation removes one vertex and creates three vertices, adds three edges, and adds one face. Hence $P_k$ has $2k$ vertices, $3k$ edges, and $k+2$ faces, verifying~\textup{(\ref{ih:counts})}.

The dihedral angles along surviving edges are unchanged. Each new edge lies in $F_k\cap F_j$ for some $j\in\{k-3,k-2,k-1\}$. By Lemma~\ref{lem: dihedral}, its dihedral angle satisfies
\[
  \angle(F_k,F_j)
    =\arccos\bigl(-\langle\vec{n}_k,\vec{n}_j\rangle\bigr)
    <\frac\pi2+\arcsin(4/R).
\]
Together with the induction hypothesis, this verifies~\textup{(\ref{ih:angles})} and completes the induction.
	\end{proof}

\bibliographystyle{plain}
\bibliography{ref.bib}

\end{document}